\documentclass[a4paper,11pt,oneside,reqno]{amsart}
\usepackage[margin=3cm]{geometry}
\usepackage{amsfonts,amssymb,amsmath,amsthm}
\usepackage{graphicx, color, psfrag,hyperref}

\usepackage[english]{babel}
\newtheorem{theorem}{Theorem}

\newtheorem{lemma}{Lemma}

\newcommand{\beq}{\begin{eqnarray}}
\newcommand{\eeq}{\end{eqnarray}}

\newcommand{\beqt}{\begin{eqnarray*}}
\newcommand{\eeqt}{\end{eqnarray*}}

\newcommand{\be}{\begin{equation}}
\newcommand{\ee}{\end{equation}}

\newcommand{\bl}{\begin{lemma}}
\newcommand{\el}{\end{lemma}}

\newcommand{\bpr}{\begin{proof}}
\newcommand{\epr}{\end{proof}}

\newcommand{\bi}{\begin{itemize}}
\newcommand{\ei}{\end{itemize}}

\newcommand{\ben}{\begin{enumerate}}
\newcommand{\een}{\end{enumerate}}

\newcommand{\Z}{\mathbb Z}
\newcommand{\R}{\mathbb R}

\newcommand{\state}[3]{\langle #1 \rangle^{#2}_{#3}}

\begin{document}

\title[Absence of Dobrushin states for $2d$ long-range Ising models]{Absence of Dobrushin states \\ for $2d$ long-range Ising models}

\author{Loren Coquille}
\address{Loren Coquille, Univ. Grenoble Alpes, CNRS, Institut Fourier, F-38000 Grenoble, France}
\email{loren.coquille@univ-grenoble-alpes.fr}

\author{Aernout C.D. van Enter}
\address{Aernout C.D. van Enter, Johann Bernoulli Institute for Mathematics and Computer Science, Nijenborgh 9, 9747AG, University of Groningen, Groningen, Netherlands}
\email{a.c.d.van.enter@rug.nl}

\author{Arnaud Le Ny}
\address{Arnaud Le Ny, LAMA UMR CNRS 8050, UPEC, Universit\'e Paris-Est, 61 Avenue du G\'en\'eral de Gaulle,  94010 Cr\'eteil cedex, France, and Eurandom, TU/e Eindhoven, Den Dolech 12, 5600 MB Eindhoven, Netherlands}
\email{arnaud.le-ny@u-pec.fr}

\author{Wioletta M. Ruszel}
\address{Wioletta M. Ruszel, Delft Institute of Applied Mathematics, Technical University Delft, Van Mourik Broekmanweg 6, 2628XE, Delft, Netherlands}
\email{w.m.ruszel@tudelft.nl}


\thanks{\textit{Acknowledgments.} LC and ALN have been partially supported by the CNRS PEPS project "Ising" and visitor grants by "Networks" and "European Women in Math". ALN  has benefited from the support of CNRS, Eurandom, and the Dutch Gravitation grant "Networks". We would like to thank Yvan Velenik for his comments and encouragements.
}

\maketitle

\begin{abstract}
	We consider the two-dimensional Ising model with long-range pair inter\-actions of the form $J_{xy}\sim|x-y|^{-\alpha}$ with $\alpha>2$, mostly when $J_{xy} \geq 0$. We show that Dobrushin states ({\em i.e.} extremal non-translation-invariant Gibbs states selected by mixed $\pm$-boundary conditions) do not exist. We discuss possible extensions of this result in the direction of the Aizenman-Higuchi theorem, or concerning fluctuations of interfaces. We also mention the existence of rigid interfaces in two long-range anisotropic contexts.
\end{abstract}


\section{Introduction}

 We are interested in the possible existence of so-called {\em interface states}  (or {\em Dobrushin states}) for long-range Ising models in dimension $d=2$. Dobrushin states are extremal infinite-volume Gibbs measures selected by mixed $\pm$-boundary conditions originally described in \cite{Dob} for the standard (nearest-neighbour) Ising model in dimension three.
 
Depending on the question one asks, a long-range interaction can  behave similarly or not to nearest-neighbour ($n.n.$) models. It is well-known that for $n.n.$\ interactions,  interface states do not exist, see e.g.\  \cite{A,Gal72,MMSole77,H,CoqVel2010, DS85}. On the other hand there exist such extremal {\em and} non-translation-invariant Gibbs states in $d\geq3$ \cite{vB,Dob}. 

{In this note we consider three different examples of long-range models, with interactions which are either {isotropic}, {long-range horizontally and $n.n.$ vertically}\ or  {bi-axial long-range} (i.e.\ long-range in both horizontal and vertical directions with possibly different decays). By long-range interaction we mean pair interactions with {\em coupling constants} given by
 $$J_{{xy}}\sim {|{x-y} |}^{- \alpha}, \;\text{ with } \alpha > 2, \;  $$ 
for all $x,y \in \Z^2$, where $|\cdot|$ denotes the Euclidean norm. 

We first prove that for all $\alpha>2$ there are no Dobrushin states in the isotropic case, see Theorem \ref{Thm:iso}.}
This result is similar to what happens for short-range models. However, the precise statements we can prove are weaker than what is known for $n.n.$\ models, in particular regarding the full picture of the convex set of Gibbs measures. The Aizenman-Higuchi theorem indeed states \cite{A,H} that all infinite-volume Gibbs measures of the $2d$ ferromagnetic $n.n.$\ Ising model are convex combinations of the pure $+$ and $-$ phases.
As we shall see, we have either a statement for a subset of boundary conditions, or at low enough temperatures, but we provide directions to study these problems in more generality.\\
Note that the case of fast decays $\alpha >3$ falls within the framework of the Gertzik-Pirogov-Sinai theory. Indeed, in this decay range, a first-moment condition on the interaction applies, and allows the proof of the existence of a spontaneous magnetisation  for Ising models with $n.n$ interaction plus a long-range perturbation \cite{GGR66, Gert}. Moreover, in such a framework, it could be shown, according to Dobrushin-Shlosman (see e.g.\ the review \cite{DS85}) that all the Gibbs measures are translation-invariant, in particular Dobrushin states do not exist. The novelty of our result thus concerns the regime $2<\alpha\leq3$.\\
Let us mention that the Dobrushin-Shlosman analysis has been extended to finite-range Kac potentials at low temperature, with small enough Kac parameter \cite{Merola}. However, in the Kac limit, $\alpha=4$ represents in some sense a critical value, because a second-moment condition on the interactions stops applying: the magnetisation of ferromagnetic vector models becomes non-zero \cite{KP}, and there can occur oscillatory phases if the long-range part of the interaction is repulsive (anti-ferromagnetic). For an example of this, see e.g. \cite{PisT}. Another example of problematic Kac limits 
can be found in \cite{GatPen1970}.

{We next prove that in some long-range anisotropic cases Dobrushin states do exist, see Theorem \ref{Thm:biax}.}
 Indeed, it is  fairly easy to see that interface states do exist in a two-dimensional model with $n.n.$\ pair interactions in one direction and {\em Dyson-like} ({\em i.e.}\ one-dimensional polynomially decaying) long-range  interactions in the other direction. In this case, the elegant proof of van~Beijeren \cite{vB,BLPO79} applies; it makes use of a duplicate set of variables due to Percus and the existence of a lower-dimensional spontaneous magnetization. 
Moreover, this phenomenon also occurs for bi-axial models, where along the horizontal or the vertical axis (or both), interactions are decaying slowly enough.

{The note is organized as follows, in Section \ref{sec:def} we give the precise definitions of the models. The main results and an outline of the proofs are given in Section \ref{sec:results}.  The following Section \ref{sec:iso} and \ref{sec:anisotropic} are dedicated to the proofs of Theorem \ref{Thm:iso} and Theorem \ref{Thm:biax} respectively.}

\section{Definition of the models}\label{sec:def}

All our results concern pair potentials with  ferromagnetic interactions for Ising models, where the configuration space is $\Omega=\{-1,+1\}^{\mathbb{Z}^2}$, and coupling constants $\big(J_{xy}\big)_{x,y \in \mathbb{Z}^2}$ satisfying the standard summability condition 
\begin{equation}\label{summability}
\forall x \in \Z^2,\; \sum_{y \in \Z^2} | J_{xy} | < \infty.
\end{equation}  We focus on ferromagnetic coupling, {\em i.e.} such that $J_{xy} \geq 0$ for all $x,y \in \Z^2$. Note  that when $J_{xy}\sim |x-y|^{-\alpha}$, then $\alpha>2$ ensures \eqref{summability}. The importance of the ferromagnetic character of the models for our results and arguments actually depends on $\alpha$ (it is in particular relevant for $\alpha \in (2,3)$). 

More precisely, we define the Hamiltonian $H_\Lambda$ in a finite-volume $\Lambda\Subset\Z^2$ with \emph{boundary condition} (b.c.) $\omega\in\Omega$ as
\be \label{Hambc}
\forall \sigma \in \Omega,\; H_\Lambda (\sigma | \omega) =\frac12 \sum_{x,y \in \Lambda} J_{xy} \sigma_x \sigma_y + \sum_{x \in \Lambda,\; y  \in \Lambda^c} J_{xy} \sigma_x \omega_y.
\ee

The associated finite-volume \emph{Gibbs measure} in $\Lambda \Subset \mathbb{Z}^2$ at  \emph{inverse temperature} $\beta>0$, with boundary condition $\omega\in\Omega$ is given by
\be
\mu^{\omega}_{\Lambda}(\sigma)=\frac1{Z^{\omega}_{\Lambda}}{e^{-\beta H_{\Lambda}(\sigma|\omega)}},
\ee
where $Z^{\omega}_{\Lambda}$ is the \emph{partition function}
\be
Z^{\omega}_{\Lambda} = \sum_{\sigma \in \Omega}e^{-\beta H_{\Lambda}(\sigma|\omega)}.
\ee

We call an infinite-volume \emph{Gibbs state} (or \emph{measure}) at inverse temperature $\beta>0$ any probability measure $\mu$ on $\{-1,+1\}^{\Z^2}$ which satisfies the {\em Dobrushin-Lanford-Ruelle equations}:
\be \label{DLR}
\forall \Lambda\subset \Z^2 \text{ such that } |\Lambda|<\infty,\quad \mu (\cdot) = \int \mu^\omega_\Lambda(\cdot )d\mu(\omega).
\ee
We call $\mathcal G(\beta)$ the set of Gibbs states at inverse temperature $\beta>0.$ And we say that the model undergoes a \emph{phase transition} if $|\mathcal{G}(\beta)|> 1$. It is known that phase transitions hold at low temperature $\beta>\beta_c(\alpha)$, with $\beta_c(\alpha)\in(0,\infty)$ for $\alpha>2$.
This can be seen in various ways, either by using Griffiths inequalities \cite{Grif67-III} in the ferromagnetic case for any $\alpha >2$, or by analytic techniques \cite{GGR66} or by an extension of the Peierls argument  \cite{Gert} for $\alpha >3$,  or by reflection-positivity \cite{FILS78} for $\alpha \in (2,4)$. 
%

 It is standard (see {\cite{FV16,Geo}) that $\mathcal{G}(\beta)$ is a convex set (actually a Choquet simplex) and that weak limits of finite-volume Gibbs measures belong to $\mathcal{G}(\beta)$. More precisely, for a sequence $(\Lambda_n)_{n\ge 1}$ of finite sets in $\mathbb{Z}^2$, we write $\Lambda_n\uparrow \mathbb{Z}^2$ if for every $x\in \mathbb{Z}^2$, there exists $n_x\ge 1$, such that $x\in \Lambda_n$ for every $n\ge n_x$. 
 We say that  $\lim_{\Lambda_n\uparrow \mathbb{Z}^2}\mu^{\omega_n}_{\Lambda_n}=\mu $ if for any local function $f:\Omega\to\R$ (i.e.\ which depends on a finite number of spins), 
 	\be
 	\lim_{\Lambda_n\uparrow \mathbb{Z}^2}\langle f \rangle^{\omega}_{\Lambda_n} = \langle f \rangle^{\omega}
 	\ee
 	where $\langle \cdot \rangle^{\omega}_{\Lambda_n}$ denotes the expectation w.r.t.\ $\mu^{\omega}_{\Lambda_n}$ and $\langle \cdot \rangle^{\omega}$ the expectation w.r.t.\ $\mu$.
 	Let $\mathcal{F}$ be the sigma-algebra over $\Omega$ generated by the cylinder sets, and let $\mathcal{M}(\Omega,\mathcal{F})$ be the set of all probability measures on $(\Omega, \mathcal{F})$.  It holds that
 	\be
 	\mathcal{G}(\beta)\supseteq\left\{ \mu\in \mathcal{M}(\Omega,\mathcal{F}): \text{ there exist }(\Lambda_n)_{n\geq 1} \text{ and }(\omega_n)_{n\geq 1} \text{ s.t. } \lim_{\Lambda_n\uparrow \mathbb{Z}}\mu^{\omega_n}_{\Lambda_n}=\mu \right\},
 	\ee
 	the question whether equality of the two sets holds is discussed in \cite{Coq2015}.
 	
 	A Gibbs state $\mu\in\mathcal{G}(\beta)$ is said to be translation invariant if for any local function $f:\Omega\to\R$, and any translation $\theta:\Z^2\to\Z^2$ (there exists $v\in\Z^2$ such that $\theta(x)=x+v$, $\forall x\in\Z^2$), we have
 	$$\mu(\theta^{-1}f)=\mu(f).$$

 	When $\omega\equiv+1$ (resp.\ $\omega\equiv-1$) we write $\mu^+_\Lambda$ (resp.\ $\mu^-_\Lambda$), and the corresponding pure phase $\mu^+$ (resp.\ $\mu^-$).
 	We use the subscript $L$, and write $\mu^\omega_L$, resp.\ $\langle \cdot \rangle^\omega_L$, for the finite-volume measures and expectations on square boxes $\Lambda=\Lambda_L=([-L,+L] \cap \Z)^2$. 
Sometimes it is   useful
to also consider rectangular boxes $\Lambda_{L,M}$, of width $L$ and height $M$, centered at the origin. In this case we write $\mu_{L,M}^\omega$, resp. $\langle \cdot \rangle_{L,M}^\omega$.

We write $\omega=(\pm,h)$ for the so-called Dobrushin b.c.\ centered at height $h\in\Z$:
\begin{align}
\omega_{(i,j)}=
\left\{
\begin{matrix}
+1, & \text{ on } \{ (i,j)\in\Z^2: j\geq h\}\\
-1, & \text{ on } \{ (i,j)\in\Z^2: j< h\}.
\end{matrix}
\right.
\end{align}
Any (sub-sequential) weak limit of sequences ${(\mu^{(\pm,h)}_{\Lambda_L})}_{L\geq1}$ is written $\mu^{(\pm,h)}$. The Gibbs state $\mu^{(\pm,h)}$ is called a Dobrushin state if it is extremal and is \emph{not} translation-invariant.

\subsection{Three Examples}

{In this subsection we precise the three different models we are considering.}

\subsubsection{Model I - isotropic long-range Ising models}\label{mod:iso}
Consider  classical $2d$ extensions of long-range Dyson models, with an isotropic pair potential so that 
$J_{xy}\sim|x-y|^{-\alpha}$.
Take for example for all $ x,y \in \mathbb{Z}^2,$
\begin{equation}\label{Jisotropic}
J_{xy} =  {|x-y|^{-\alpha}},\; 
\end{equation}
where $|x-y|=\sqrt{  |x_1 - y_1|^2 + |x_2- y_2|^2 }$.
We call \emph{Model $I.1$} the case  $\alpha >3$ and  
\emph{Model $I.2$}  the case  $2 < \alpha \leq  3$.

\subsubsection{Model II - anisotropic long-range/$n.n.$ Ising models}\label{mod:aniso}
Consider a mixed long-range and $n.n.$ translation-invariant interaction whose interactions  are  $n.n.$ vertically and 'Dyson-like' horizontally, {\em i.e.} 
of the form 
\begin{align}
J_{xy} &= 1 \quad \text{ if }  x_1=y_1  \text{ and } y_2=x_2 \pm 1,\nonumber \\
J_{xy} &= 0 \quad \text{ if } |y_2-x_2| > 1, \text{ or } y_2 = x_2 \pm 1 \text{ and } x_1 \neq y_1,\nonumber\\
J_{xy} &= |x_1-y_1|^{-\alpha_1} \quad  \text{ if }  x_2=y_2.
\end{align}

\subsubsection{Model III - bi-axial, possibly anisotropic, long-range Ising models}\label{mod:biax}
Consider  'Dyson-like' long-range interactions in both horizontal and vertical directions with possibly different decays  $\alpha_1,\alpha_2>1$ and $\alpha_1\in(1,2)$: 
\begin{align}
J_{xy} &=|x_2-y_2|^{-\alpha_2} \quad  \text{ if }   x_1=y_1,\nonumber\\
J_{xy} &= |x_1-y_1|^{-\alpha_1} \quad  \text{ if }  x_2=y_2,\nonumber\\
J_{xy}&=0 \quad \text{ otherwise}.
\end{align}

\section{Main results}\label{sec:results}

Let us state our main results and discuss the ideas of the proofs. 

\begin{theorem}\label{Thm:iso}
In the case of Model I (isotropic long-range) defined in \ref{mod:iso}, for any $\alpha>2$ and any $\beta<\infty$, Dobrushin states do not exist. Furthermore for $\alpha>3$ and $\beta$ large enough, all Gibbs measures are translation-invariant.
\end{theorem}

\begin{theorem}{(Corollaries of \cite{vB} and \cite[Appendix B]{BLPO79})}\label{Thm:biax}\emph{}\\
	\vspace{-.5cm}
\begin{enumerate}
\item In the case of Model II (anisotropic long-range/$n.n.$) defined in \ref{mod:aniso}, there exist Dobrushin states for any $\alpha_1 \in (1,2)$ provided $\beta>\beta_c(\alpha_1,d=1)$.
\item In the case of Model III (bi-axial long-range) defined in \ref{mod:biax}, there exist Dobrushin states for any $\alpha_1 \in (1,2)$ and $\alpha_2>1$ provided $\beta>\beta_c(\alpha_1,d=1)$.
\end{enumerate}
\end{theorem}

Our approach to prove Theorem \ref{Thm:iso} distinguishes between a $n.n.$-like picture for $\alpha>3$, and a more long-range one in the case $2 < \alpha \leq  3$. In the first case, our work is in accordance with the previous analyses contained in works of Ginibre {\em et al.}\ \cite{GGR66}, Gertzik \cite{Gert}, Dobrushin-Shlosman \cite{DS85} or Bricmont-Lebowitz-Pfister \cite{BLP79}.  In the second case, we use arguments close to those
mainly introduced  by Fr\"ohlich-Pfister for continuous spin systems \cite{FrPf81, FrPf86}, also described in the review of Dobrushin-Shlosman \cite{DS85}.

Theorem \ref{Thm:biax}(1) can be proved using van Beijeren's original proof \cite{vB}, while point (2) is a particular case of the models treated in the Appendix B of \cite{BLPO79}; their arguments imply the existence of Dobrushin states states either in one or two directions, depending on whether one or two of the decay rates $\alpha_1, \alpha_2$ is between 1 and 2. These results can be extended to long-range models, as long as some symmetries are kept, see Section \ref{sec:anisotropic}.

\section{Proof of Theorem \ref{Thm:iso}: Absence of Dobrushin states in the isotropic case}\label{sec:iso}

We consider the \emph{Model I}, defined in \ref{mod:iso}, in the phase transition region (i.e.\ at inverse temperature $\beta > \beta_c(\alpha,d=2)$). 
 
 \subsection{Model I.2 : Fast decays, $\alpha > 3$}\emph{}\\
The main observation we use is that the difference between two Dobrushin b.c.\ located at different horizontal heights is obtained by flipping all spins in two half-lines. If the maximal energy between a half-line left of the origin and a half-plane right of the origin, to which it is perpendicular,  is uniformly bounded, the arguments of \cite{BLP79} apply and we can conclude that there is no pure interface Gibbs state. We remark as an aside that this argument does not need the ferromagnetic character of the model.

\subsubsection{Energy difference in a finite box with Dobrushin boundary conditions}\emph{}\\
In the box $\Lambda=\Lambda_L$, we write 
\begin{align}\label{up-down}
\Lambda^{up,h}=\{ (i,j): j\geq h\} \cap \Lambda^c, \quad 
\Lambda^{down,h}=\{ (i,j): j< h\} \cap \Lambda^c.
\end{align}
We have
\[
- H^{\pm,0}_{\Lambda}(\sigma) = \frac12\sum_{x,y\in \Lambda} \sigma_x \sigma_y J_{xy}+ \sum_{x\in \Lambda, y\in \Lambda^{up,0}}\sigma_x  J_{xy} - \sum_{x\in \Lambda, y\in \Lambda^{down,0}}\sigma_x  J_{xy}
,\]
\[
- {H}^{\pm,1}_{\Lambda}(\sigma) = \frac12\sum_{x,y\in \Lambda} \sigma_x \sigma_y J_{xy}+ \sum_{x\in \Lambda, y\in \Lambda^{up,1}}\sigma_x  J_{xy} - \sum_{x\in \Lambda, y\in \Lambda^{down,1}}\sigma_x  J_{xy}
,\]
Thus,
{
\begin{multline*} 
	|H^{\pm,0}_{\Lambda}(\sigma) - {H}^{\pm,1}_{\Lambda}(\sigma)| 
 =\\
 \left |\sum_{x\in \Lambda, y\in \Lambda^{up,0}}\sigma_x  J_{xy}  
 - \sum_{x\in \Lambda, y\in \Lambda^{down,0}}\sigma_x  J_{xy}
 - \sum_{x\in \Lambda, y\in \Lambda^{up,1}}\sigma_x  J_{xy}
 +  \sum_{x\in \Lambda, y\in \Lambda^{down,1}}\sigma_x  J_{xy}  \right | ,
 \end{multline*}}
so  by writing $x=(i_x,j_x), y=(i_y,j_y)$, and using the explicit expression \eqref{Jisotropic} of  $J_{xy}$ we get
 \begin{eqnarray*}
|H^{\pm,0}_{\Lambda}(\sigma)- {H}^{\pm,1}_{\Lambda}(\sigma)|  & \leq  &
\sum_{(i_y,0) \in \Lambda^c}  \sum_{(i_x,j_x)\in \Lambda} |(i_x-i_y)^2+j_x^2 |^{-\alpha/2}\\
& =&  \sum_{(i_y,0) \in \Lambda^c}  \sum_{i_x=-L}^L \sum_{j_x=-L}^L |(i_x-i_y)^2+j_x^2 |^{-\alpha/2} \\
& \leq &  \sum_{(i_y,0) \in \Lambda^c}  \sum_{i_x=-L}^L \sum_{j_x=-\infty}^{\infty} |(i_x-i_y)^2+j_x^2 |^{-\alpha/2}\\
& \leq & C \sum_{(i_y,0) \in \Lambda^c}  \sum_{i_x=-L}^L  |i_x-i_y|^{1-\alpha} \\
& \leq & C \sum_{i_y=L+1}^{\infty}  \sum_{i_x=0}^L \left ( (i_y-i_x)^{1-\alpha} + (i_x+i_y)^{1-\alpha}\right ) ,
\end{eqnarray*}
 with possibly different constants $C$ from line to line. Now provided $\alpha >3$ we have
\[
\begin{split}
\sum_{i_y=L+1}^{\infty} \sum_{i_x=0}^L (i_y-i_x)^{1-\alpha} <\infty,
\end{split}
\]
and analogously
\[
\sum_{i_y=L+1}^{\infty} \sum_{i_x=0}^L (i_y+i_x)^{1-\alpha} <\infty.
\]
It follows that
\begin{equation}\label{TranslateEnergy}
|H^{\pm,0}_{\Lambda}(\sigma)- {H}^{\pm,1}_{\Lambda}(\sigma)| = O(1),
\end{equation}
 which is uniformly bounded.

Notice, as indicated before,  that the estimates here do not depend on the sign of the interactions. Thus they also work for 
decaying interactions
of any sign.
\smallskip

For slower decays $2 < \alpha \leq 3$, which will be studied in Section \ref{subsec-alpha23} , this estimate becomes unbounded  as $L$ grows. It is moreover sharp in the sense that the energy difference  between a plus configuration and a minus configuration on the  half-line 
$\{x  < 0, y=0 \}$ interacting with the plus configuration on  the half-plane 
$\{ x \geq 0 \}$ indeed is infinite.

\subsubsection{Absence of Dobrushin States}\emph{}\\
Once we have the energy estimate (\ref{TranslateEnergy}), the proof goes as in \cite{BLP79}:
 finite energy difference implies that the states obtained as a weak limit with $(\pm,0)$ and $(\pm,1)$ Dobrushin b.c.\ 
 are absolutely continuous  with respect to each other. Therefore they have the same components in their extremal decomposition and they are equal if one of them is extremal. So if the limit $\lim_{\Lambda\uparrow\Z^2}\mu_{\Lambda}^{(\pm,0)}$ is extremal, then it coincides with $\lim_{\Lambda\uparrow\Z^2}\mu_{\Lambda}^{(\pm,1)}$ and  is thus translation-invariant : Dobrushin states are excluded for any $\alpha>3$ and any $\beta<\infty$.
 
\subsubsection{Translation-Invariance at Low Temperature}\emph{}\\
 The case of fast decays $\alpha >3$ falls within the framework of the Gertzik-Pirogov-Sinai theory. Indeed, those models satisfy a Peierls condition at low enough temperature as shown in \cite{Gert}. In such a framework, all the Gibbs measures should be translation-invariant, as claimed and described in the review \cite{DS85}. From this, coupled with the fact recently extended to more general contexts  by Raoufi  \cite{Raoufi17} that the $\mu^+$ and $\mu^-$ states are the only translation-invariant extremal states, one gets also the convex decompositions in terms of these pure states.
For example, $$\mu^{(\pm,0)}=\lim_{\Lambda\uparrow\Z^2}\mu_{\Lambda}^{(\pm,0)}=\frac{1}{2} (\mu^- +  \mu^+).$$
 


\subsection{Model I.2 : Slow decays, $2< \alpha \leq 3$}\label{subsec-alpha23}\emph{}\\
In this case, although the maximal interaction energy between a half-line left of the origin and a half-plane right of it is infinite, we show that the expected interaction energy in a state with Dobrushin boundary conditions still remains finite. We use here both the ``anti-symmetry'' between upper and lower parts of the box and the ferromagnetic character of the interaction. Indeed, the argument breaks down if the interaction has for example   alternating signs in the vertical direction. 
\subsubsection{Energy Difference between the Dobrushin ground-state and the Dobrushin ground state flipped on a half-line}\emph{}\\
Split the lattice $\mathbb{Z}^2$ into 
\begin{align}
A^+&=\{ (i,j): j\geq 1 \}\cup \{ (i,0): i> 0\},\nonumber\\
A^-&=\{ (i,j): j\leq -1 \},\nonumber\\
A^0&=\{ (i,0): i\leq0 \}. 
\end{align}

We define $\sigma_{GS}$ to be the ground state of the $(\pm,0)$ Dobrushin boundary condition, that is the configuration consisting of $+1$ in $A^+\cup A^0$ and $-1$ in $A^-$.
We call $\sigma_{GS,step}$ the configuration $\sigma_{GS}$ which is flipped on the half line $A^0$, that is
consisting in $+1$ in $A^+$ and $-1$ in $A^0\cup A^-$. Then
\[
-H(\sigma_{GS}) = \frac12\sum_{x,y\in A^+} J_{xy} + \frac12\sum_{x,y \in A^-}J_{xy} - \sum_{x\in A^+, y\in A^-} J_{xy} + \sum_{x\in A^+, y\in A^0}  J_{xy} -\sum_{x\in A^0, y\in A^-} J_{xy}
\]
and
\[
-H(\sigma_{GS,step}) = \frac12\sum_{x,y\in A^+} J_{xy} + \frac12\sum_{x,y \in A^-} J_{xy} - \sum_{x\in A^+, y\in A^-} J_{xy} \textcolor{red}{-} \sum_{x\in A^+, y\in A^0}  J_{xy} \textcolor{red}{+}\sum_{x\in A^0, y\in A^-} J_{xy}
\]
writing as before $x=(i_x,j_x), y=(i_y,j_y)$, the energy difference is equal to
\[
\begin{split}
\big|H(\sigma_{GS}) - H(\sigma_{GS,step})\big|& = 2 \left| \sum_{x\in A^+, y\in A^0}  J_{xy} -  \sum_{y \in A^0, x\in A^-} J_{xy}\right| \\
& = \left | \sum_{i_y = - \infty}^0 \sum_{i_x\in \mathbb{Z}} \sum_{j_x=1}^{\infty} J_{xy} + \sum_{i_y=-\infty}^0 \sum_{i_x=1}^{\infty} J_{xy} - \sum_{i_y = - \infty}^0 \sum_{i_x\in \mathbb{Z}} \sum_{j_x=-\infty}^{-1} J_{xy}\right |.
\end{split}
\]
By symmetry of the couplings $J_{xy}$, the first and third term cancel each other out. Thus, provided $\alpha>2$,
\[
\begin{split}
\big|H(\sigma_{GS}) - H(\sigma_{GS,step})\big|& =  
\sum_{i_y=-\infty}^0 \sum_{i_x=1}^{\infty} J_{xy} =  
\sum_{i_y=-\infty}^0 \sum_{i_x=1}^{\infty} |i_x-i_y|^{-\alpha} \\
& = \sum_{i_y =0}^{\infty} \sum_{i_x=1}^{\infty} (i_x+i_y)^{-\alpha}<\infty.
\end{split}
\]

In words, the argument uses the fact that the interaction of the negative 
half-line 
$\{ i < 0, j=0 \}$ and the positive half-line $\{ i \geq 0, j=0 \}$ is finite, while the interaction of the half-line with any plus spin above the line is canceled by the interaction with the reflected minus spin below the line.

\subsubsection{ Positive Temperatures : Absence of Dobrushin states}\emph{}\\
A similar argument will still hold for expected energy differences at positive temperatures.
Indeed, the interaction energy of a spin at distance $\ell$ from a half-plane, interacting with it,  is maximally $O(\ell^{2- \alpha}$), but its expectation in the Gibbs state with Dobrushin b.c. at more or less the same height is $O(\ell^{1- \alpha})$.
Summing over the line just above the Dobrushin interface at level 1/2 gives that the total expected energy cost of shifting is uniformly bounded, thus the relative entropy
  between the two Dobrushin states is finite, and thus again,  they are the same once they are extremal \cite{BLP79}. We use here  a general strategy of Pfister, inspired by previous work of Araki, more precisely formalised in \cite{FrPf86} that yields estimates and results at any positive temperature.

More precisely, let us compute the relative entropy between measures with Dobrushin b.c. $(\pm,0)$ and $(\pm,1)$, which is the expectation of the energy difference computed in one of the states, see e.g. \cite{FrPf81}.
It is thus given by
\begin{align}
\Delta {\mathcal H}(\mu^{\pm,0},\mu^{\pm,1})
&:=\lim_{L\to\infty}\lim_{M\to\infty}\langle|H^{\pm,0}_\Lambda(\sigma)-H^{\pm,1}_\Lambda(\sigma)|\rangle^{\pm,1}_{L,M}\nonumber\\
&=\lim_{L\to\infty}\langle|H^{\pm,0}_{S_L}(\sigma)-H^{\pm,1}_{S_L}(\sigma)|\rangle^{\pm,1}_{L,\infty}
\end{align}
where $S_L=\{(i,j):-L\leq i\leq L,j\in\mathbb Z\}$ is the vertical strip of width $2L+1$.
Note that the measure  $\mu_{L,\infty}^{\pm,1}$ obtained as a weak limit of the measure $\mu_{L,M}^{\pm,1}$ as  the height of the box $M \to \infty$,  has the following (anti-symmetry) property of expectations:
\begin{equation}\label{antisym}
\langle \sigma_{(i,j)}\rangle^{\pm,1}_{L,\infty}= - \langle \sigma_{(i,1-j)}\rangle^{\pm,1}_{L,\infty}.
\end{equation}

Let us define $A^{left}_0=\{ (i,0): i\leq -L \}$ and $A^{right}_0=\{ (i,0): i\geq L \}$ and write as before $x=(i_x,j_x), y=(i_y,j_y)$. 

We bound
\begin{align}
& \langle|H^{\pm,0}_S(\sigma) -H^{\pm,1}_S(\sigma)|\rangle^{\pm,1}_{L,\infty}\nonumber\\
&=\sum_{x\in S}\sum_{y \in A_0^{left}\cup A_0^{right}} J_{xy} \cdot \langle \sigma_x\rangle^{\pm,1}_{L,\infty}
=2\sum_{x\in S}\sum_{y \in A_0^{right}} J_{xy} \cdot \langle \sigma_x\rangle^{\pm,1}_{L,\infty}\nonumber\\
&=2\sum_{i_x=-L}^L\sum_{j_x=-\infty}^\infty\sum_{i_y=L}^\infty J_{xy} \cdot\langle \sigma_x\rangle^{\pm,1}_{L,\infty}\nonumber\\
&=2\sum_{i_x=-L}^L\sum_{j_x=1}^\infty\sum_{i_y=L}^\infty \left(J_{(i_x,j_x),(i_y,0)} \cdot\langle \sigma_{(i_x,j_x)}\rangle^{\pm,1}_{L,\infty}+J_{(i_x,1-j_x),(i_y,0)} \cdot\langle \sigma_{(i_x,1-j_x)}\rangle^{\pm,1}_{L,\infty}\right)\nonumber\\
&\label{line}= 2\sum_{i_x=-L}^L\sum_{j_x=1}^\infty\sum_{i_y=L}^\infty \left(J_{(i_x,j_x),(i_y,0)} -J_{(i_x,1-j_x),(i_y,0)} \right) \langle \sigma_{(i_x,j_x)}\rangle^{\pm,1}_{L,\infty}\\
& \leq 2 
	\sum_{i_x=-L}^L\sum_{j_x=1}^\infty\sum_{i_y=L}^{L+\ell} J_{xy}+\nonumber\\
 &\hspace{3cm}+ 2\sum_{i_x=-L}^L\sum_{i_y=L+\ell}^{\infty}
\sum_{j_x=1}^\infty
 \left({{((i_x-i_y)^2+j_x^2)}^{-\frac \alpha 2}} -{((i_x-i_y)^2+(1-j_x)^2)^{-\frac \alpha 2} } 
\right)\nonumber\\
&<\infty \label{estim}
\end{align}

where \eqref{line} comes from the anti-symmetry property \eqref{antisym}, and \eqref{estim} is valid for $\ell$ large enough (vertically symmetric spins are almost at the same distance to $y$ when $\ell$ is large, so the last term is integrable in $j_x$ and only two integrations are left).\\

Hence the relative entropy remains  uniformly bounded as $L$ grows to infinity for any $\alpha>2$ and any $\beta<\infty$. This implies the absolute continuity of the weak limits got by these b.c.'s so eventually translation-invariance of any extremal measures obtained in this manner, and in particular the absence of any Dobrushin states : the weak limits we consider cannot be both extremal and non-translation-invariant.
\subsubsection{About Translation-Invariance}\emph{}\\
Concerning translation-invariance itself, as our models do not fall anymore within the range of decays satisfying the Gertzik-Pirogov-Sinai condition, we cannot use low temperature Dobrushin-Shlosman results. 


To conclude that the measures resulting from  convex decompositions of extremal states  are translation-invariant requires to exclude possibly other non-translation-invariant extremal states. However,  we believe that our energy estimate, or a variation thereof, will provide translation-invariant  measures in the thermodynamic limit for all boundary conditions having a finite number of sign changes. Note that rigid diagonal interfaces are expected to exist for nearest-neighbour models  in dimension four \cite{MMSole77}.
To exclude other extremal non-translation-invariant states,
we suspect that we need extra arguments which in the nearest-neighbour case were provided  by Aizenman \cite{A} or Higuchi \cite{H} after percolation results of Russo \cite{Russo78} and techniques based on correlation inequalities developed in the seventies (see \cite{GMSole72, MMSole75, MMSole77, Leb77, H1} and references therein). To get the full convex picture in the  long-range case in dimension two, we probably need to investigate the validity of Russo's results (which roughly reduces the problem of excluding non-translation-invariance along axes) and afterwards the validity of its extension by Higuchi \cite{H}, which probably requires new ideas.

\section{Proof of Theorem \ref{Thm:biax}: Dobrushin states for anisotropic models}\label{sec:anisotropic}

Our rigidity results for anisotropic two-dimensional long-range models will be
applications of the results
of Bricmont {\em et al.} \cite{BLPO79}, so we give  it here and sketch the ingredients of the proof (which is an extension of van Beijeren's one \cite{vB}).

Consider   the general Ising model with coupling constants $\big( J_{xy} \big)_{x,y \in {\mathbb{Z}^d}}$ being
\begin{enumerate}
\item{{\it Ferromagnetic :}} $J_{xy} \geq 0$.
\item{{\it Reflection-invariant :}}  For any $x=(x_1,\dots,x_d)$, let $\bar{x}=(-x_1,\dots,x_d)$, then
$$
J_{xy}=J_{\bar{x}\bar{y}}.
$$
\item{{\it Growing in the quarter-plane :}} $J_{xy} \geq J_{x\bar{y}}$ $\forall x_1, y_1\geq 0$.
\item{{\it Summable }} For any $x \in \mathbb{Z}^d$,
\begin{equation*}\label{BLOPvB}
\sum_{y \in \mathbb{Z}^d} J_{xy} < \infty.
\end{equation*}
\end{enumerate}

Note that our Models II and III (defined in \ref{mod:aniso} and \ref{mod:biax} respectively) fulfill these conditions. To get rigidity of the interface created by the Dobrushin b.c.\ $(\pm,0)$, the idea of van Beijeren is to lower-bound the magnetization (in the full box) of a spin lying at height 1, by its magnetization in a $(d-1)$-dimensional  box with $+$ b.c. 
The first Hamiltonian writes (keeping the notations \eqref{up-down})
$$-H(\sigma)= \frac12\sum_{x,y \in \Lambda} J_{xy} \sigma_x \sigma_y + \sum_{x\in \Lambda, y\in\Lambda^{up,1}} J_{xy} \sigma_x - \sum_{x\in \Lambda, y\in\Lambda^{down,1}} J_{xy} \sigma_x,$$
while the latter corresponds to an Ising model on the plane $\Lambda_0=\{x\in\Z^d:x_1=0\}$ with the same infinite-range interaction within it, 
$$
-H'(\sigma')= \frac12\sum_{x,y \in \Lambda_0} J_{xy} \sigma'_x \sigma'_y + \sum_{x \in \Lambda_0, y \in \Lambda^c, y_1=0} J_{xy} \sigma'_x.
$$
Then one uses a clever change of variables, originally due to Percus, which allows to involve notably the differences $t_x=\sigma_x-\sigma_x'$ (with $x\in\Lambda_0$), and express the sum of the two Hamiltonians $H(\sigma)+H'(\sigma')$ as a ferromagnetic Hamiltonian in the new variables. As a consequence, GKS correlation inequalities apply, and we get the following result (cf.\  Appendix B of \cite{BLPO79}):
\begin{align*}
\state{t_x}{\pm,+}{\Lambda,\Lambda_0}&\geq0\quad \forall x \in \Lambda_0,\;  
\end{align*}
hence,
\begin{align*}
\state{\sigma_x}{\pm}{\Lambda} &\geq \state{\sigma'_x}{+}{\Lambda_0}.
\end{align*}

The expectation of the $t_x$ variables is taken under the measure associated to the sum of the two Hamiltonians, while the expectation of $\sigma_x'$ is performed under the $(d-1)$-dimensional Gibbs states with $+$ b.c. at the same temperature. Thus, as soon as spontaneous magnetization occurs for the latter,  this implies a strictly positive magnetization of $\sigma_x$ and thus the existence of a non-translation-invariant Gibbs state in dimension $d$. 
Notice that this lower-dimensional phase transition condition is not fulfilled in the isotropic long-range models treated above, because their well-definedness requires $\alpha >2$, for which there is no phase transition in dimension one. To get such a phase transition and positive magnetization, one has to consider very long ranges in dimension one with decays $1< \alpha_1 \leq 2$.
This has motivated the introduction of our anisotropic Models II and III, for which this lower-dimensional spontaneous magnetization holds at low temperature. 
%

\section{Conclusion}

We have given some conditions under which two-dimensional interface Gibbs states exist, and also some conditions under which a plausible recipe of making them does not work. This falls short of the statement, which is known for nearest-neigbour models, that all Gibbs states for the Ising model are translation-invariant.  

Also, for nearest-neighbour models it is known how big the interface fluctuations are, that the existence of a phase transition is equivalent to the existence of a positive surface tension, and that the existence of interface states is equivalent to the positivity of a step free energy.
Which of these statements, possibly in a modified form,  hold for long-range models is not clear to us at the moment. 

What about the size of the interface fluctuations?
Our results imply that the fluctuations of the interface cannot stay uniformly bounded. For $\alpha > 3$, at very low temperature, it is to be expected that they show a diffusive behavior  of order $\sqrt L$, like other models in the Peierls regime. The situation for $2< \alpha \leq 3$ seems more delicate, as was suggested by Fr\"ohlich and Zegarlinski \cite{FZ91}. The authors indeed predict a behavior of order 
$O(L^{\alpha-2 \over 2})$ 
in a one-dimensional discrete Gaussian model with long-range interactions which seems to be a reasonable effective interface model in the low-temperature case. Their argument looks  somewhat similar to an upper bound on correlations (which implies a lower bound on fluctuations) developed in \cite{COE,Ent87}, for long-range spin-glasses, respectively one-dimensional long-range spin models. However, whereas those papers provided a rigorous but non-sharp bound, the arguments of \cite{FZ91} at this point appear to be potentially sharp, but, as the authors say, there still needs some hard work, even in the presumably
simpler discrete Gaussian model, to make them rigorous.


\bibliographystyle{abbrv}
\bibliography{library-long-range}

\end{document}